\documentclass[11pt]{article}
\usepackage{graphicx} 

\usepackage{fancyhdr}
\usepackage{isomath}
\usepackage{mathtools} 
\usepackage{amsbsy}
\usepackage{amssymb}
\usepackage{amscd,amsfonts}
\usepackage{bigints}
\usepackage{graphicx}
\usepackage{verbatim}
\usepackage{euscript}
\usepackage{alltt}
\usepackage{stmaryrd}
\usepackage{relsize}
\usepackage{enumerate}
\usepackage{url}
\usepackage[stable]{footmisc}
\usepackage{comment}
\usepackage[font=small,labelfont=bf]{caption}
\usepackage[font=small,labelfont=bf]{subcaption}
\usepackage{float}
\usepackage{mwe}
\usepackage{algorithm}
\usepackage{algpseudocode}
\usepackage{xcolor}
\usepackage[super]{nth}
\usepackage{soul}

\DeclareGraphicsExtensions{.eps,.pdf}
\usepackage{amsmath}
\usepackage{amsbsy}
\usepackage{amssymb}
\usepackage{amscd}
\usepackage{amsfonts}

\newcommand{\R}{\mathbb R}

\newcommand{\beq}{\begin{equation}}
\newcommand{\eeq}{\end{equation}}
\newcommand{\beqs}{\begin{eqnarray}}
\newcommand{\eeqs}{\end{eqnarray}}
\newcommand{\beql}{\begin{equation} \label}
\newcommand{\half}{\frac{1}{2}}

\newcommand{\calA}{{\cal A}}

\newcommand{\calO}{{\cal O}}

\newcommand{\veps}{\varepsilon}

\newcommand{\p}{\partial}

\usepackage[margin=1in]{geometry}

\newcommand{\dee}{\mathcal{D}}
\newcommand{\scl}{\mathcal{L}}

\newcommand{\M}{\mathbb M}
\newcommand{\B}{\mathbb B}
\newcommand{\K}{\mathbb K}
\newcommand{\G}{\Gamma}

\date{}
\begin{document}
\title{Ideal Magnetohydrodynamics and Field Dislocation Mechanics}

\author{Amit Acharya\thanks{Department of Civil \& Environmental Engineering, and Center for Nonlinear Analysis, Carnegie Mellon University, Pittsburgh, PA 15213, email: acharyaamit@cmu.edu.}}

\maketitle

$\quad $ \textit{This paper is dedicated to Professor Luc Tartar on the occasion of his $75^{th}$ birthday.\\}
\begin{abstract}
\noindent The fully nonlinear (geometric and material) system of Field Dislocation Mechanics is reviewed to establish an exact analogy with the equations of ideal magnetohydrodynamics (ideal MHD) under suitable physically simplifying circumstances. Weak solutions with various conservation properties have been established for ideal MHD recently by Faraco, Lindberg, and Sz\'ekelyhidi using the techniques of compensated compactness of Tartar and Murat and convex integration; by the established analogy, these results would seem to be transferable to the idealization of Field Dislocation Mechanics considered. A dual variational principle is designed and discussed for this system of PDE, with the technique transferable to the study of MHD as well.

\end{abstract}

\section{Introduction}
Dislocations are line defects in the perfect crystal structure of many solids. They are found in abundance in most common metals, and exist either as closed loops or curves that cannot end within the body, often forming complex tangles as they interact and change shape through the stresses they produce, due to the disturbance produced by them in the perfect crystal structure, and due to applied loads. These elastic interactions are of long-range, e.g.~in 2-d given by $\sim r^{-1}$ (away from their cores), where $r$ is their separation. It is the motion of large numbers of dislocations that produces large plastic, or permanent, deformation in metals, this ductility being of significant technological importance. At the same time, their tangling inhibits this motion giving the material its strength. Low ductility makes for less resilient brittle materials, prone to catastrophic failures by cracking even in the presence of high strength. High-strength, high-ductility materials is a holy grail of material design.

Dislocations have been studied as topological defects in classical elasticity, initiated by Volterra in 1907 \cite{volterra1907equilibre}, with a great deal understood about their elastic stresses and static, energetic interactions in linear elasticity \cite{hirth_lothe}. Their dynamic aspects, in particular the connection to the macroscopic plasticity of metals, remains a difficult subject and a matter of current study. A complete dynamical theory of unrestricted material and geometric nonlinearity within
continuum mechanics  has taken shape in recent years (see, e.g., \cite{action_3, acharya2019structure, zhang2015single, arora2020finite, arora2020unification,arora_acharya_ijss} with detailed bibliographic threads to earlier works in the overall subject). Within this particular framework, rigorous mathematical studies have been undertaken in \cite{acharya2010travelling, acharya2011equation, ach_slemrod}. In this paper, an analogy with ideal magnetohydrodynamics is established in Sec.~\ref{sec:analogy}, and some ideas for defining variational dual solutions to the PDE of dislocation mechanics are discussed in Sec.~\ref{sec:vp}.

\section{The parallel between ideal MHD and FDM}\label{sec:analogy}
We will always employ the summation convention on repeated indices, unless explicitly mentioned otherwise. A Rectangular Cartesian coordinate system is employed in all instances, and $X$ represents the alternating tensor with components represented as $e_{ijk}$.

Let $\Omega \in \R^3$ be a fixed domain and $[0,T]$ an interval of time. For the fields
\begin{equation*}
    \begin{aligned}
        \rho &: \Omega \times [0,T] \to \R^+, \qquad & \mbox{mass density}\\
        v &: \Omega \times [0,T] \to \R^3, \qquad & \mbox{material velocity}\\
        V &: \Omega \times [0,T] \to \R^3, \qquad & \mbox{dislocation velocity}\\
        T &: \Omega \times [0,T] \to \R^{3 \times 3}_{sym}, \qquad & \mbox{Cauchy stress}\\
        \alpha &: \Omega \times [0,T] \to \R^{3 \times 3}, \qquad & \mbox{Dislocation density}\\
        W &: \Omega \times [0,T] \to \R^{3 \times 3}, \qquad & \mbox{inverse elastic distortion},
    \end{aligned}
\end{equation*}
the equations of Field Dislocation Mechanics are given by
\begin{subequations}
    \begin{align}
        & \p_t(\rho v) + div (\rho v \otimes v)  = div \, T; &\qquad \p_t (\rho v_i) + \p_j (\rho v_i v_j) &- \p_j T_{ij}  = 0 \notag\\
        & & &\mbox{(balance of linear momentum)} \\
        & \p_t \alpha + curl(\alpha \times(v+V)) = 0;  &\qquad   \p_t \alpha_{ij} + e_{jrk} & \p_r( e_{kmn} \alpha_{im} (v_n + V_n)) = 0   \notag\\
        & & &\mbox{(Burgers vector conservation)} \label{eq:disloc_evol}\\
        & curl \, W  = - \alpha; & \qquad  e_{ijk} \p_j W_{rk}  = - \alpha_{ri}  \qquad &\mbox{(incompatibility of elastic distortion)} \label{eq:incomp}\\
        & \p_t \rho + div(\rho v) = 0; & \qquad  \p_t \rho + \p_i (\rho v_i) = 0 \qquad & \mbox{(balance of mass)}
    \end{align}
\end{subequations}
along with the constitutive equations (for free energy density given by $\phi(W) + \veps |\alpha|^2$, a sum of elastic and core energy densities)
\begin{equation}
    \begin{aligned}
        T & = - \rho W^T \p_W \phi - \veps \rho |\alpha|^2 I + \veps \rho \alpha^T \alpha; \qquad &  T_{ij} = -  \rho W_{ki} \p_{W_{ij}} \phi - \veps \rho (\alpha_{rs}\alpha_{rs}) \delta_{ij} + \veps \rho \alpha_{ki}\alpha_{kj}\\
        V & = M X\left( \left(- \rho \p_W \phi + \veps \rho curl \alpha \right)^T \alpha \right); \qquad & V_i = M_{ij} e_{jmn} ( (- \rho \, \p_{W_{rm}} \phi + \veps \rho \, e_{mpq} \p_p \alpha_{rq}) \alpha_{rn} ),
    \end{aligned}
\end{equation}
where $M$ is a constant, positive semi-definite, dislocation mobility tensor, $\veps \ll 1$ is a constant scalar, and the `Burgers vector' is the topological charge of a dislocation line. The statements \eqref{eq:disloc_evol} and \eqref{eq:incomp}, along with the physical assumption that plastic strain rate arises only from the presence of dislocation motion implies the statement that the elastic distortion evolves according to
\[
\p_t W + \nabla W v + W \nabla v - \alpha \times V = 0; \qquad \p_t W_{ij} + v_k \p_k W_{ij} + W_{im} \p_j v_m - (\alpha \times V)_{ij} = 0,
\]
which is an alternate and, strictly speaking, more general, evolution statement in lieu of \eqref{eq:disloc_evol} in FDM.

We now invoke a set of assumptions to probe the behavior of FDM in a restricted setting. These are:
\begin{itemize}
    \item The time scales of interest are much faster than the characteristic time scales of dislocation motion relative to the material built into the dislocation mobility tensor $M$ defining $V$, so that one may assume that we  are in the regime of elastic, i.e.~non-dissipative, \textit{ideal} dislocation dynamics where $\alpha  \times V  \approx 0$ (it can be shown that the dissipation of mechanical energy in the model is characterized by the density $T:(\alpha \times V)$). This also implies that the initial dislocation density profile is transported solely by the material velocity during any deformation process being studied;
    
    \item the deviatoric stress response is a bounded function of elastic distortion, whereas the hydrostatic stress response is extremely stiff; to simplify matters we replace this stiff hydrostatic stress response with the constraint of flow incompressibility, $\p_i v_i = 0$, and introduce a pressure field $p$ as the corresponding Lagrange multiplier to be solved for as an independent field. We then satisfy balance of mass by choosing $\rho = \rho_0$, a constant;
    
    \item We assume that the inverse elastic distortion field has rapid spatial variation of the form $W(x,t) \sim \tilde{W}\left(\veps^{-\left(\frac{1}{2} + m\right)} x,t \right), m > 0$ (inherited, in part, from the initial dislocation density field), with the bounded deviatoric stress response being $O(1)$, so that the latter can be ignored in writing the stress as $T \approx - p I + \alpha^T\alpha$ (after defining $\alpha := \veps^{-m} \sqrt{\rho_0} \, curl \, \tilde{W}$ and absorbing all hydrostatic parts of the stress into the constitutively undetermined pressure field).
\end{itemize}
With this understanding, we have, for
\[
v: \Omega \times [0,T] \to \R^3; \quad \alpha: \Omega \times [0,T] \to \R^{3 \times 3}; \quad p: \Omega \times [0,T] \to \R; \quad B:\Omega \times [0,T] \to \R^3,
\]
where $B$ is the external magnetic field in MHD, the correspondence shown in Table \ref{tab:parallel} between ideal FDM and ideal MHD (with $(v,p)$ being distinct fields between FDM and MHD):
\begin{table}[h]
\def\arraystretch{1.3}
    \centering
    \begin{tabular}{|c|c|}
    \hline
       Ideal FDM  &  Ideal MHD\\
       \hline
        $\p_t v + div \left( v \otimes v - \alpha^T \alpha + p I \right) = 0 \qquad$ & $ \qquad \p_t v + div ( v \otimes v - B \otimes B + p I) = 0$\\
        \hline
        $div \, v = 0 \qquad$ & $ \qquad div \, v = 0$\\
        \hline
        $\p_t \alpha + curl (\alpha \times v) = 0 \qquad$ & $\qquad \p_t B + curl (B \times v) = 0$\\
        \hline
        $div \, \alpha = 0 \qquad$ & $ \qquad div \, B = 0$\\
        \hline
    \end{tabular}
    \caption{Correspondence between ideal FDM and MHD}
    \label{tab:parallel}
\end{table}

For a concise description of ideal MHD, see, e.g., \cite{faraco2022rigorous}. We note here that for $v$ specified, the dislocation transport equation in FDM constitutes three uncoupled sets of vectorial equations for the rows of the matrix $\alpha$ (the $curl$ on a matrix argument there amounts to row-wise $curl$s of the matrix). This does not constitute a substantial difference with ideal MHD, with the primary difference arising from the nonlinearities $B_i B_j$ and $\alpha_{ki} \alpha_{kj}$ in balance of linear momentum.

Defining the field $\chi: \Omega \times [0,T] \to \R^{3 \times 3}$ in FDM such that
\[
curl \, \chi = \alpha
\]
one can define the analog of magnetic helicity in FDM as the quantity
\[
\int_\Omega \chi : \alpha \, dx = \int_\Omega \chi_{ij} \alpha_{ij} \, dx.
\]
As discussed by Faraco, Lindberg, and Sz\'ekelyhidi \cite{faraco2022rigorous}, it is shown by them that magnetic helicity, $\int_\Omega A\cdot B \, dx$, is  conserved for ideal limits of resistive, viscous MHD, if such limits exist (for a definition of the ideal limit of resistive, viscous MHD, see \cite{faraco2022rigorous}). They also show that magnetic helicity can be conserved with energy dissipation for some $L^\infty$ solutions $(u, B)$ to ideal MHD, as well as the fact that for $ 2\leq p < 3$ there are fields $(u,B) \in (0,T; L^p(\mathbb{T}^3))$ that do not conserve magnetic helicity or energy.  Given the above parallel between ideal MHD and ideal FDM, one might expect similar results to exist for the latter as well. These authors also discuss the even more difficult question of proving the existence of ideal limit solutions of resistive, viscous MHD, and appropriate admissibility criteria on ideal MHD that can allow the recovery of only ideal limit solutions of resistive, viscous MHD. In view of such results and questions, the correspondence established above already allows for many of the mathematical results of ideal MHD to be potentially transferrable to ideal FDM. In the following, a class of variational principles is proposed for the study of the PDE system of ideal FDM, with similar formalism also being applicable to the equations of ideal MHD. The approach \cite{action_2,action_3} bears strong links to the ideas of `hidden convexity' for PDE advanced by Y.~Brenier \cite{brenier_book,brenier2018initial}, with a primary philosophical difference being our reliance on, and exploitation of, the Euler-Lagrange equations of the variational principle, with the deliberate recognition that this is not an impediment of any sort.

\section{Variational principle for ideal FDM}\label{sec:vp}
We consider the system of PDE of ideal FDM on a domain $\Omega \in \R^3$ given by
\begin{subequations}\label{eq:ifdm}
    \begin{align}
        \p_t v_i + \p_j (v_i v_j - \alpha_{ki}\alpha_{kj} + p \delta_{ij}) & = 0 \notag\\
        \p_t \alpha_{ij} + e_{jrp}  \p_r( e_{pms} \alpha_{im} v_s) & = 0 \label{eq:alpha_dot}\\
        \p_i v_i & = 0, \notag
    \end{align}
\end{subequations}
noting that \eqref{eq:alpha_dot} implies $div\, \alpha = 0$. Since the question of appropriate boundary conditions for ideal FDM is non-trivial, we will assume the domain $\Omega = (0,1) \times (0,1) \times (0,1)$ to be the unit cube with periodic boundary condition on the fluxes, i.e. with no sum implied on the index $j$,
\begin{equation}\label{eq:primal_bc}
    \begin{aligned}
        (v_i v_j - \alpha_{ki}\alpha_{kj} + p \delta_{ij})\big|_{(x_j = 1,t)} - (v_i v_j - \alpha_{ki}\alpha_{kj} + p \delta_{ij})\big|_{(x_j = 0,t)} & = 0 \\
        (e_{pjr}e_{pms} \alpha_{im}v_s)\big|_{(x_j = 1,t)} - (e_{pjr}e_{pms} \alpha_{im}v_s)\big|_{(x_j = 0,t)} & = 0.
    \end{aligned}
\end{equation}
The initial conditions specified on $\Omega$ are
\begin{equation}\label{eq:primal_ic}
    v_i(x,0) = v_i^{(0)}(x); \qquad \alpha_{ij}(x,0) = \alpha_{ij}^{(0)}(x) \qquad x \in \Omega.
\end{equation}
We refer to these equations and side-conditions as the `primal' system and adopt the perhaps somewhat unconventional point of view of treating these equations as constraints for optimizing a more-or-less arbitrarily chosen \textit{auxiliary potential} $H(U,x,t)$ of the primal variables and space-time, with a requirement defined next. 

Let
\begin{equation*}
    \begin{aligned}
        & \calO := \Omega \times (0,T); \qquad (v, \lambda, \bar{v}): \calO \to \R^3; \qquad (\alpha, A, \bar{\alpha}): \calO \to \R^{3 \times 3}; \qquad (p,\mu, \bar{p}):\calO \to \R \\
        & U:=(v,\alpha, p) \quad \mbox{primal fields}; \qquad D := (\lambda, A, \mu) \quad \mbox{dual fields}; \qquad \bar{U}:=(\bar{v},\bar{\alpha}, \bar{p}) \quad \mbox{base states}\\
        & \dee := (\p_t \lambda, \nabla \lambda, \nabla \mu, \p_t A, \nabla A).
    \end{aligned}
\end{equation*}
Defining a \textit{Lagrangian} as
\begin{equation}\label{eq:Lagrangian}
\begin{aligned}
    \scl_H(U,\dee,x,t) &: = - v_i \p_t \lambda_i - (v_i v_j - \alpha_{ki}\alpha_{kj} + p \delta_{ij})\p_j \lambda_i\\ 
   & \quad \   - v_i \p_i \mu - \alpha_{ij} \p_t A_{ij} - e_{pjr}e_{pms} \alpha_{im} v_s \p_r A_{ij} + H(U, x,t)
\end{aligned}
\end{equation}
(obtained, up to addition of $H$, by forming the scalar products of \eqref{eq:ifdm} with the dual (Lagrange multiplier) fields and integrating by parts), the requirement on the potential $H$ is that $\p_U \scl_H (U,D,x,t)= 0$ be `solvable' for $U$ in terms of $(\dee, x,t)$. I.e., the choice of $H$ must be such as to enable the existence of a function $U^{(H)}(\dee, x,t)$ such that
\begin{equation}\label{eq:dldu}
\frac{\p \scl_H}{\p U} \left( U^{(H)}(\dee, x,t), \dee, x,t\right) = 0
\end{equation}
is satisfied for $(x,t) \in \calO$ and for each $\dee = \dee(x,t)$ arising from an appropriately non-trivial class of dual fields of interest.

We refer to the function $U^{(H)}$ as a \textit{dual-to-primal} (DtP) mapping.

When the DtP mapping can be defined, further defining a `dual' functional corresponding to the PDE system \eqref{eq:ifdm} as
\begin{equation}
    \label{eq:dual_ifdm}
    S_H[D] = \int_\calO \scl_H \left(U^{(H)}(\dee(x,t),x,t),\dee(x,t), x,t \right) \, dx dt  - \int_\Omega \lambda_i(x,0) v_i^{(0)}(x) \,dx - \int_\Omega A_{ij}(x,0) \alpha_{ij}^{(0)}(x) \, dx
\end{equation}
with the conditions
\begin{subequations}\label{eq:dual_ifdm_dirichlet}
\begin{align}
    & \lambda\big|_{(x_j = 1,t)} = \lambda\big|_{(x_j = 0,t)}; \qquad \alpha\big|_{(x_j = 1,t)} = \alpha\big|_{(x_j = 0,t)}, \qquad x \in \p \Omega, \ j = 1,2,3; \notag\\
    & \mu\big|_{\p \Omega \times (0,T)} \ \mbox{specified arbitrarily};  \notag\\
    & \lambda\big|_{\Omega \times \{T\}}, A\big|_{\Omega \times \{T\}} \ \mbox{specified arbitrarily}, \label{eq:final_time_bc}
\end{align}  
\end{subequations}
it can be checked that, \textit{due to \eqref{eq:dldu} and the necessarily affine dependence of $\scl_H$ on $\dee$, the Euler-Lagrange (E-L) equations and natural boundary conditions of the dual functional $S_H[D]$ are the primal statements \eqref{eq:ifdm}-\eqref{eq:primal_bc}-\eqref{eq:primal_ic} with the replacement}
\begin{equation}\label{eq:fundamental}
    U \to \hat{U}^{(H)}, \qquad \mbox{where} \qquad \hat{U}^{(H)}(x,t) : = U^{(H)}(\dee(x,t), x,t).
\end{equation}
Clearly, the purpose of the auxiliary potential $H$ is to facilitate a \textit{change of variables} $U^{(H)}$ by dominating any non-monotonicity in the variable $U$ arising from the primal PDE in $\p_U \scl_H$ by a sufficiently monotone $\p_U H$, achieved through a choice of $H$ with a sufficiently strongly positive-definite Hessian. Of course, the matter is not entirely as clear-cut in that $\scl_H$ depends on $\dee$ which $H$ cannot depend upon, but to the extent that in the course of solving any specific problem, the structure of the PDE and the choice of $H$ can be made to work together to ensure the existence of a such a mapping (at least for a sufficiently large class of dual fields), one can obtain a weak variational formulation of ideal FDM with some rather desirable properties, as we discuss next.

Given the quadratic nonlinearity of \eqref{eq:ifdm}, we now choose to work with a potential $H$ of the form:
\begin{equation}\label{eq:H}
H(U,x,t) := \frac{1}{2}\left( a_v |v - \bar{v}(x,t)|^2 + a_\alpha |\alpha - \bar{\alpha}(x,t)|^2 + a_p |p - \bar{p}(x,t)|^2 \right),
\end{equation}
where $ 1 \ll a_v, a_\alpha, a_p \in \R^+$ are arbitrary constants, and the base state $\bar{U}$ comprises freely assignable functions over space-time (and hence `states,' even though time-dependent). Then the system $\frac{\p \scl}{\p U} = 0$ (we drop the subscript $H$ now, having made a specific choice) becomes an algebraic system to be solved at each $(x,t) \in \calO$ for $(\hat{v}, \hat{\alpha}, \hat{p})|_{(x,t)}$, given $(\dee(x,t), \bar{U}(x,t))$:
\begin{equation}\label{eq:mapping}
    \begin{aligned}
        \p_{v_i} \scl & = \quad a_v (\hat{v}_i - \bar{v}_i) - (\p_j \lambda_i + \p_i \lambda_j) \hat{v}_j \quad - e_{pjr}e_{pmi} \hat{\alpha}_{km} \p_r A_{kj} \quad  - \p_t \lambda_ i - \p_i \mu  = 0\\
        \p_{\alpha_{ij}} \scl & = \quad - e_{pkr} e_{pjs}\hat{v}_s \p_r A_{ik}  \quad + a_\alpha (\hat{\alpha}_{ij} - \bar{\alpha}_{ij}) + (\p_s \lambda_j + \p_j \lambda_s) \hat{\alpha}_{is} \quad - \p_t A_{ij}  = 0\\
        \p_p \scl & =  \quad a_p(\hat{p} - \bar{p})  \quad  - \p_i \lambda_i  = 0.
    \end{aligned}
\end{equation}
This system is to be satisfied simultaneously with the E-L equation of the dual functional, \eqref{eq:ifdm} with the replacement $U \to \hat{U}$ which we now refer to as $\hat{\eqref{eq:ifdm}}$. One typical iteration of a conceptual algorithm for the nonlinear coupled problem (that has been successfully implemented \cite{KA1, a_arora_thesis, sga, KA2}) consists of making a guess for the fields $D$, evaluating $\hat{U}$ from the mapping \eqref{eq:mapping}, and solving the linearization of $\hat{\eqref{eq:ifdm}}$ to improve the guess for $D$ until convergence in the solution of $\hat{\eqref{eq:ifdm}}$, i.e.~a standard Newton method for which the starting guess $D_{guess} = 0$ suffices with corresponding $\hat{U}_{guess} = \bar{U}$ satisfying the DtP mapping.

We develop next the explicit form of the dual functional for this problem - which also provides a guideline for a general PDE system with at most quadratic nonlinearities. Indexing the elements of $U$ 
in a one-dimensional array by upper-case Latin letters and similarly the elements of $\dee$ 
by upper-case Greek letters (the dependence of $\dee$ on $\nabla \lambda$ is only through its symmetric part and that on $\nabla A$ is through its $curl$),
the Lagrangian in \eqref{eq:Lagrangian} can be expressed, in terms of constant matrices $\M_{I\Gamma}$ and $\B_{\Gamma JK}$ with $\B$ symmetric in last two indices, as
\begin{equation}\label{eq:expl_Lagr_1}
    \begin{aligned}
        \scl & = U \cdot \M \dee + \half \dee \cdot \B : (U \otimes U) + \half (U - \bar{U}) \cdot a \, (U - \bar{U})\\
        & = U_I \M_{I \Gamma} \dee_\Gamma + \half \dee_\Gamma \B_{\Gamma JK} U_J U_K  + \half a_{IJ}(U_I - \bar{U}_I)(U_J - \bar{U}_J)\\
        & = (U - \bar{U}_I)\M_{I\Gamma} \dee_\Gamma  +  \half a_{IJ}(U_I - \bar{U}_I)(U_J - \bar{U}_J) + \half \dee_\G \B_{\G IK}(U_I - \bar{U}_I)(U_K - \bar{U}_K)\\
        & \quad + \bar{U}_I \M_{I \Gamma} \dee_\G + \half \dee_\G \B_{\G JK} U_J \bar{U}_K + \half \dee_\G \B_{\G JK} \bar{U}_J U_K - \half \dee_\G \B_{\G JK} \bar{U}_J \bar{U}_K.
    \end{aligned}
\end{equation}
Then the DtP mapping is to be generated from
\begin{equation}\label{eq:expl_Lagr_2}
\begin{aligned}
    \p_U \scl & = \M \dee + a (U - \bar{U}) + \dee \cdot \B (U - \bar{U}) + \dee \cdot \B \bar{U} = 0\\
    \p_{U_I} \scl & = \M_{I \G} \dee_\G + a_{IJ}(U_J - \bar{U}_J) + \dee_\G \B_{\G IK}(U_K - \bar{U}_K) + \dee_\G \B_{\G IK} \bar{U}_K = 0.
\end{aligned} 
\end{equation}
Contracting \eqref{eq:expl_Lagr_2} with $(U_I - \bar{U}_I)$ and substituting for $(U_I - \bar{U}_I) \M_{I \G} \dee_\G$ in \eqref{eq:expl_Lagr_1} one obtains
\begin{equation*}
    \begin{aligned}
        \scl & = - (U_I - \bar{U}_I) a_{IJ} (U_J - \bar{U}_J) + \half a_{IJ}(U_I - \bar{U}_I)(U_J - \bar{U}_J)\\
        & \quad - \dee_\G \B_{\G IK}(U_I - \bar{U}_I)(U_K - \bar{U}_K) + \half \dee_\G \B_{\G IK}(U_I - \bar{U}_I)(U_K - \bar{U}_K)\\
        & \quad - (U_I - \bar{U}_I) \dee_\G \B_{\G IK} \bar{U}_K + U_I \dee_\G \B_{\G IK}  \bar{U}_K\\
        & \quad + \bar{U}_I \M_{I \Gamma} \dee_\G - \half \dee_\G \B_{\G JK} \bar{U}_J \bar{U}_K\\
        & = - \half (U_I - \bar{U}_I) a_{IJ} (U_J - \bar{U}_J) - \half \dee_\G \B_{\G IK}(U_I - \bar{U}_I)(U_K - \bar{U}_K)\\
        & \quad + \bar{U}_I \M_{I \Gamma} \dee_\G + \half \dee_\G \B_{\G JK} \bar{U}_J \bar{U}_K\\
        & = - \half (U - \bar{U}) \cdot a \, (U - \bar{U}) - \half \dee \cdot \B : \left(\left(U - \bar{U} \right) \otimes \left(U - \bar{U} \right) \right) + \bar{U} \cdot \M \dee + \half \dee \cdot \B : (\bar{U} \otimes \bar{U}).
    \end{aligned}
\end{equation*}
Using the DtP mapping once again in the form
\begin{equation}\label{eq:K}
    \K_{IJ}\big|_{\dee} (U_J - \bar{U}_J) = - (\M_{I \G} + \B_{\G I K} \bar{U}_K) \dee_\G \ ; \qquad \K_{IJ}\big|_{\dee} := a_{IJ} + \dee_\G \B_{\G IJ},
\end{equation}
one obtains, with the definition
\[
\calA^{\Pi \G}\big|_{(\dee,\bar{U})} := \left(\M_{J \Pi} + \B_{\Pi J K} \bar{U}_K \right)  \left( \K\big|_\dee^{-1} \right)_{JI} \, \left(\M_{I \G} + \B_{\G I K} \bar{U}_K \right),
\]
\begin{equation*}
    \begin{aligned}
        \scl & = - \half  \dee_\Pi \cdot \calA^{\Pi \G}\big|_{(\dee,\bar{U})} \dee_\G  + \bar{U}_I \M_{I \Gamma} \dee_\G + \half \dee_\G \B_{\G JK} \bar{U}_J \bar{U}_K\\
        & =: - \half \dee \cdot \calA\big|_{(\dee,\bar{U})} \dee + \bar{U} \cdot \M \dee + \half \dee \cdot \B : (\bar{U} \otimes  \bar{U}).
    \end{aligned}
\end{equation*}

Correspondingly, the dual functional for ideal FDM may be explicitly written as
\begin{equation}\label{eq:explicit_dual_func}
\begin{aligned}
    S[D] & = \int_\calO  - \half \dee \cdot \calA\big|_{(\dee,\bar{U})} \dee + \bar{U} \cdot \M \dee + \half \dee \cdot \B : (\bar{U} \otimes  \bar{U}) \, dx dt  \\
    & \quad - \int_\Omega \lambda(x,0) \cdot v^{(0)}(x) \,dx - \int_\Omega A(x,0) : \alpha^{(0)}(x) \, dx,\\
    & \qquad \mbox{along with the conditions \eqref{eq:dual_ifdm_dirichlet}.}
\end{aligned}
\end{equation}

We make the following observations:
\begin{itemize}
    \item One may wonder why it should be possible to solve a primal initial value problem as a boundary value problem in the time-like direction? Indeed, all classical experience suggests that any such endeavor most likely produces an ill-posed problem, see, e.g., in this specific context of solving a dual variational problem the remarks of Brenier \cite[p.~39]{brenier_book} ``Anyway, it seems foolish to solve the IVP
problem by a space-time convex minimization technique. Indeed, ..... we are very likely to get optimality equations of space-time elliptic type and therefore ill posed," as well as in \cite[p.~580]{brenier2018initial}. The point to note here is that the boundary value problem in the time-like direction, corresponding to a primal initial-value-problem arises \textit{in the dual fields} which satisfy second-order degenerate elliptic equations that admit final-time conditions - this has been established in \cite[Sec.~7]{action_2}-\cite{KA1} in the context of linear and nonlinear equations. A simple explanation of this in the context of the present work is as follows: we note that imposition of the boundary conditions in time \eqref{eq:final_time_bc} allows the mapped primal fields $(\hat{u}, \hat{\alpha})$ to take on the values dictated by the primal initial value problem by virtue of the fact that the mapping equations \eqref{eq:mapping} depend on $(\p_t \lambda, \p_t A)$ and specifying $(\lambda(x, T), A(x,T)), x \in \Omega$, does not restrict the values of $(\p_t \lambda(x,T), \p_t A(x,T)), x \in \Omega$.
    \item We note that if the base state $\bar{U}$ in \eqref{eq:H} is a weak solution of \eqref{eq:ifdm}, then \eqref{eq:mapping} and the fundamental property of the scheme \eqref{eq:fundamental} show that the dual state $D = 0$ is a critical point of the dual functional. Thus, corresponding to every weak solution of the primal system at least one dual functional can be designed which has a critical point that generates a global-in-time solution to the primal problem through the corresponding DtP mapping. Moroever, since the dual functional is concave, this solution is also a global mazimizer. This may be thought of as a consistency check for the scheme.
    \item When the primal problem has uniqueness of solutions, if a solution exists to the dual problem for any choice of the auxiliary potential $H$, then such a dual solution defines the unique primal solution through the DtP mapping. When the primal system has non-unique solutions, then the choice of $H$ is expected to act like a `selection criterion.' Following the arguments provided in \cite[Sec.~3]{dual_cont_mech_plas} it can be shown that the dual PDE system $\hat{\eqref{eq:ifdm}}$ is locally degenerate elliptic in a neighborhood around the dual state $D = 0$. Thus, it may be reasonable to expect that if a solution to the primal system happens to be `close' to the base state $\bar{U}$ that corresponds to the $D = 0$ dual state, then such a primal solution may be recoverable by working with a shifted quadratic \eqref{eq:H} designed with this base state.
    \item If solving the mapping equations and subsequently using these mappings to define the dual functional were to be reinterpreted as performing an $\inf_U \scl(U,\dee,x,t)$ to define the bulk integrand of \eqref{eq:explicit_dual_func}, then by the affineness of $\scl$ in $\dee$, the dual integrand is concave in $\nabla D$ regardless of the monotonicity of the left-hand-side of the primal system in the argument $U$. Hence, if the dual functional can now be shown to be coercive in the dual fields, then the existence of a maximizer can be expected. This can be thought of as a definition of a variational dual solution to the primal system, because of the formal consistency of our scheme in defining a solution to the primal PDE system through the DtP change of variables.

    At the level of the explicit dual functional \eqref{eq:explicit_dual_func}, a computation of the second variation about the state $\bar{D} = 0$ shows that the second variation is positive semi-definite. If now $\bar{U}$ corresponds to a solution of the primal problem, then the $\bar{D}=0$ state is a critical point of the dual variational problem and hence is neutrally stable. Thus, regardless of the stability of $\bar{U}$ as a solution to the primal system, it can be approached as a stable solution through the proposed dual scheme. This seems to have implications for computing unstable primal solutions through our scheme.
    \item Another intriguing possibility seems to be the question of using the base state to be a solution to the primal system with a small regularization for which solutions may be accessible and then see what solution to the ideal system (FDM or MHD) may be possible to define or approximate through our scheme. A concrete example of this ideas is worked out in \cite{KA2} in the context of computing entropy solutions to the inviscid Burgers equation by using a base state of Burgers equation with small viscosity (a curious aside in this regard is that despite similarities in ideal MHD and ideal FDM, the formal structure of their energy dissipative counterpart systems are very different in nature \cite{acharya2011equation}).
\end{itemize}

\section*{Acknowledgment}
I thank L\'aszl\'o Sz\'ekelyhidi and Vladimir \v{S}ver\'ak for helpful discussion. This work was supported by the Simons Pivot Fellowship grant \# 983171, and the Center for Extreme Events in Structurally Evolving Materials, Army Research Laboratory Contract No. W911NF2320073.

\bibliographystyle{alpha}\bibliography{ref}

\begin{thebibliography}{ZAWB15}

\bibitem[AA19]{arora_acharya_ijss}
Rajat Arora and Amit Acharya.
\newblock Dislocation pattern formation in finite deformation crystal plasticity.
\newblock {\em International Journal of Solids and Structures}, 184(2):114--135, 2020, electronically published Feb. 26, 2019.

\bibitem[AA20]{arora2020unification}
Rajat Arora and Amit Acharya.
\newblock A unification of finite deformation ${J}_2$ {V}on-{M}ises plasticity and quantitative dislocation mechanics.
\newblock {\em Journal of the Mechanics and Physics of Solids}, 143:104050, 2020.

\bibitem[Ach22]{action_2}
Amit Acharya.
\newblock Variational principles for nonlinear {PDE} systems via duality.
\newblock {\em Quarterly of Applied Mathematics}, LXXXI:127--140, 2023, Article electronically published on September 26, 2022.

\bibitem[Ach23a]{dual_cont_mech_plas}
Amit Acharya.
\newblock A hidden convexity in continuum mechanics, with application to classical, continuous-time, rate-(in)dependent plasticity.
\newblock \url{https://arxiv.org/abs/2310.03201},2023.

\bibitem[Ach23b]{action_3}
Amit Acharya.
\newblock A dual variational principle for nonlinear dislocation dynamics.
\newblock {\em Journal of Elasticity}, \url{https://doi.org/10.1007/s10659-023-09998-5}, 2023.

\bibitem[AKS19]{acharya2019structure}
Amit Acharya, Robin~J. Knops, and Jeyabal Sivaloganathan.
\newblock On the structure of linear dislocation field theory.
\newblock {\em Journal of the Mechanics and Physics of Solids}, 130:216--244, 2019.

\bibitem[AMZ10]{acharya2010travelling}
Amit Acharya, Karsten Matthies, and Johannes Zimmer.
\newblock Travelling wave solutions for a quasilinear model of field dislocation mechanics.
\newblock {\em Journal of the Mechanics and Physics of Solids}, 58(12):2043--2053, 2010.

\bibitem[Aro23]{a_arora_thesis}
Abhishek Arora.
\newblock {\em A study of nonlinear deformations and defects in the actuation of soft membranes, rupture dynamics, and mesoscale plasticity}.
\newblock Phd thesis, Carnegie Mellon University, September 2023.
\newblock Available at \url{https://www.proquest.com/docview/2869411693?accountid=9902}.

\bibitem[AS23]{ach_slemrod}
Amit Acharya and Marshall Slemrod.
\newblock Existence, uniqueness, and long-time behavior of linearized field dislocation dynamics.
\newblock {\em Quarterly of Applied Mathematics}, LXXXI:247--258, 2023.

\bibitem[AT11]{acharya2011equation}
Amit Acharya and Luc Tartar.
\newblock On an equation from the theory of field dislocation mechanics.
\newblock {\em Bollettino dell’Unione Matematica Italiana}, 9:409--444, 2011.

\bibitem[AZA20]{arora2020finite}
Rajat Arora, Xiaohan Zhang, and Amit Acharya.
\newblock Finite element approximation of finite deformation dislocation mechanics.
\newblock {\em Computer Methods in Applied Mechanics and Engineering}, 367:113076, 2020.

\bibitem[Bre18]{brenier2018initial}
Y.~Brenier.
\newblock The initial value problem for the {E}uler equations of incompressible fluids viewed as a concave maximization problem.
\newblock {\em Communications in Mathematical Physics}, 364(2):579--605, 2018.

\bibitem[Bre20]{brenier_book}
Y.~Brenier.
\newblock Examples of hidden convexity in nonlinear {PDE}s.
\newblock \url{https://hal.science/hal-02928398/document}, 2020.

\bibitem[FLS22]{faraco2022rigorous}
D.~Faraco, S.~Lindberg, and L{\'a}szl{\'o} Sz{\'e}kelyhidi.
\newblock Rigorous results on conserved and dissipated quantities in ideal {MHD} turbulence.
\newblock {\em Geophysical \& Astrophysical Fluid Dynamics}, 116(4):237--260, 2022.

\bibitem[HL82]{hirth_lothe}
J.~P. Hirth and J.~Lothe.
\newblock {\em Theory of dislocations}.
\newblock Krieger, 1982.

\bibitem[KA23]{KA1}
Uditnarayan Kouskiya and Amit Acharya.
\newblock Hidden convexity in the heat, linear transport, and {E}uler’s rigid body equations: A computational approach.
\newblock {\em to appear in Quarterly of Applied Mathematics}, \url{https://arxiv.org/abs/2304.09418}, 2023.

\bibitem[KA24]{KA2}
Uditnarayan Kouskiya and Amit Acharya.
\newblock {I}nviscid {B}urgers as a degenerate elliptic problem.
\newblock to appear in Quarterly of Applied Mathematics, \url{https://arxiv.org/abs/2401.08814}, 2024.

\bibitem[SGA24]{sga}
Siddharth Singh, Janusz Ginster, and Amit Acharya.
\newblock {A} {H}idden {C}onvexity of {N}onlinear {E}lasticity.
\newblock \url{https://arxiv.org/abs/2401.08538}, 2024.

\bibitem[Vol07]{volterra1907equilibre}
V.~Volterra.
\newblock Sur l'{\'e}quilibre des corps {\'e}lastiques multiplement connexes.
\newblock In {\em Annales scientifiques de l'{\'E}cole normale sup{\'e}rieure}, volume~24, pages 401--517, 1907.

\bibitem[ZAWB15]{zhang2015single}
Xiaohan Zhang, Amit Acharya, Noel~J. Walkington, and Jacobo Bielak.
\newblock A single theory for some quasi-static, supersonic, atomic, and tectonic scale applications of dislocations.
\newblock {\em Journal of the Mechanics and Physics of Solids}, 84:145--195, 2015.

\end{thebibliography}
\end{document}